\documentclass[11pt,reqno]{amsart}
\usepackage{amssymb,amscd,amsbsy}
\usepackage{amsthm}
\newcommand{\lb}{\linebreak}

\renewcommand{\a}{\alpha}
\renewcommand{\b}{\beta}

\renewcommand{\d}{\delta}
\newcommand{\e}{\varepsilon}

\newcommand{\z}{\zeta}

\newcommand{\s}{\sigma}
\renewcommand{\t}{\tau}

\newcommand{\f}{\varphi}
\renewcommand{\o}{\omega}

\newcommand{\D}{\Delta}
\renewcommand{\L}{\Lambda}

\renewcommand{\O}{\Omega}

\newcommand{\E}{{\mathcal E}}

\newcommand{\F}{{\mathcal F}}

\newcommand{\h}{{\mathcal H}}

\newcommand{\p}{{\mathcal P}}

\newcommand{\X}{{\mathcal X}}
\newcommand{\Y}{{\mathcal Y}}

\newcommand{\T}{{\Bbb T}}
\newcommand{\pp}{{\Bbb P}}

\newcommand{\R}{{\Bbb R}}
\newcommand{\Z}{{\Bbb Z}}

\newcommand{\0}{{\boldsymbol{0}}}

\newcommand{\m}{{\boldsymbol m}}
\newcommand{\bS}{{\boldsymbol S}}

\newcommand{\rf}[1]{(\ref{#1})}

\newcommand{\df}{\stackrel{\mathrm{def}}{=}}

\newcommand{\supp}{\operatorname{supp}}

\newcommand{\trace}{\operatorname{trace}}

\newcommand{\const}{\operatorname{const}}

\newcommand{\eeq}{\end{equation}}
\newcommand{\beq}{\begin{equation}}
\newcommand{\bay}{\begin{eqnarray}}
\newcommand{\ba}{\begin{align*}}
\newcommand{\ea}{\end{align*}}
\newcommand{\ey}{\end{eqnarray}}
\newcommand{\bey}{\begin{eqnarray*}}
\newcommand{\eey}{\end{eqnarray*}}

\newcommand{\be}{\infty}

\newcommand{\bl}{\blacksquare}

\newcommand{\Pf}{{\bf Proof. }}

\newcommand{\ov}{\overline}

\newtheorem{thm}{\hspace{\parindent}Theorem}[section]

\newtheorem{lem}[thm]{\hspace{\parindent}Lemma}

\theoremstyle{remark}

\newtheorem*{rem*}{Remark}

\newcommand\ite{{\Li\,\hat\odot_{\rm i}\,L^\be}}
\newcommand\te{{\Li\,\hat\odot\,L^\be}}

\newcommand\Li{{\rm Lip}}

\begin{document}

\newcommand{\vse}{\vspace{.2in}}
\numberwithin{equation}{section}

\title{\bf An extension of the Koplienko--Neidhardt\\ trace formulae}
\author{V.V. Peller}
\thanks{The author is partially supported by NSF grant DMS 0200712}
\maketitle

\newcommand{\mt}{{\mathcal T}}

\begin{abstract}
Koplienko \cite{Ko} found a trace formula for perturbations of
self-adjoint operators by operators of Hilbert Schmidt class $\bS_2$.
A similar formula in the case of unitary operators was obtained by Neidhardt \cite{N}.
In this paper we improve their results and obtain sharp conditions under which the 
Koplienko--Neidhardt trace formulae hold.
\end{abstract}

\setcounter{section}{0}
\section{\bf Introduction}
\setcounter{equation}{0}

\

The spectral shift function for a trace class perturbation of a self-adjoint (unitary) operator
plays a very important role in perturbation theory. It was introduced in a special case by
I.M. Lifshitz \cite{L} and in the general case by M.G. Krein \cite{Kr}.
He showed that for a pair of self-adjoint (not necessarily bounded) operators $A$ and $B$ 
satisfying $B-A\in\bS_1$ there exists a unique function $\xi\in L^1(\R)$ such that
\bay
\label{sa}
\trace\big(\f(B)-\f(A)\big)=\int_\R\f'(x)\xi(x)\,dt,
\ey
whenever $\f$ is a function on $\R$ with Fourier transform of $\f'$ in $L^1(\R)$. The function $\xi$ is
called the {\it spectral shift function corresponding to the pair} $(A,B)$. We use the notation $\bS_1$ 
for the class of nuclear operators (trace class) on Hilbert space. 

A similar result was obtained in \cite{Kr2} for pairs of unitary operators $(U,V)$ with 
\lb$V-U\in\bS_1$. For each such pair there exists a function $\xi$ on the unit circle $\T$
of class $L^1(\T)$ such that 
\bay
\label{u}
\trace\big(\f(V)-\f(U)\big)=\int_\T\f'(\z)\xi(\z)\,d\m(\z),
\ey
whenever $\f'$ has absolutely convergent Fourier series. (Throughout the paper $\m$ is normalized
Lebesgue measure on $\T$.) Such a function $\xi$ is unique modulo a constant
and it is called a {\it spectral shift function corresponding to the pair} $(U,V)$. We refer the reader
to the lectures of M.G. Krein \cite{Kr3}, in which the above results were discussed in detail
(see also the survey article \cite{BY}).

Note that spectral shift function plays an important role in perturbation theory. We mention here
the paper \cite{BK}, in which the following important formula was found:
$$
\det S(x)=e^{-2\pi{\rm i}\xi(x)},
$$
where $S$ is the scattering matrix corresponding to the pair $(A,B)$; see also the monograph \cite{Y}.

It was shown later
in \cite{BS3} that formulae \rf{sa} and \rf{u} hold under less restrictive assumptions on $\f$.

Note that the right-hand sides of \rf{sa} and \rf{u} make sense for an arbitrary Lipschitz function $\f$.
However, it turns out that the condition $\f\in\Li$ (i.e., $\f$ is a Lipschitz function) does not imply
that $\f(B)-\f(A)$ or $\f(V)-\f(U)$ belong to $\bS_1$. This is not even true for bounded $A$ and $B$ and
continuously differentiable $\f$. The first such examples were given in \cite{F}. 

In \cite{Pe1} and \cite{Pe2} with the help of the nuclearity criterion for Hankel operators (see
recent monograph \cite{P}) necessary conditions (in terms of Besov classes and Carleson
measures) were found on $\f$ for the operator $\f(B)-\f(A)$ (or
$\f(V)-\f(U)$) to belong to $\bS_1$. Those necessary conditions also imply that the condition $\f\in C^1$
is not sufficient for those operators to be in $\bS_1$ (even for bounded $A$ and $B$).

In the same papers \cite{Pe1} and \cite{Pe2} sharp sufficient conditions were found. It was shown in
\cite{Pe1} that if $\f$ is a function on $\T$ of Besov class $B_{\be1}^1$, then trace formula
\rf{u} holds. Similarly, it was shown in \cite{Pe2} that if $\f$ is a function on $\R$ of Besov class
$B_{\be1}^1(\R)$, then trace formula \rf{sa} holds. The definition of the above Besov classes will be
given in \S 2. Note that though these sufficient conditions are not necessary, the gap between the
necessary conditions and the sufficient conditions obtained in \cite{Pe1} and \cite{Pe2} is rather
narrow. Note also that in \cite{ABF} a better sufficient condition was found; however, it seems to me
that the condition $\f\in B_{\be1}^1$ is easier to work with.

In Koplienko's paper \cite{Ko} the author considered the case of perturbations of Hilbert-Schmidt class
$\bS_2$. Let $A$ and $B$ be a self-adjoint operators such that $K\df B-A\in\bS_2$. In this case the
operator $\f(B)-\f(A)$ does not have to be in $\bS_1$ even for very nice functions $\f$. The idea of
Koplienko was to consider the operator
$$
\f(B)-\f(A)-\frac{d}{ds}\Big(\f(A+sK)\Big)\Big|_{s=0}
$$
and find a trace formula under certain assumptions on $\f$. It was shown in \cite{Ko} that there exists
a unique function $\eta\in L^1(\R)$ such that
\bay
\label{ktf}
\trace\left(\f(B)-\f(A)-\frac{d}{ds}\Big(\f(A+sK)\Big)\Big|_{s=0}\right)=\int_\R\f''(x)\eta(x)\,dx
\ey
for rational functions $\f$ with poles off $\R$. The function $\eta$ is called the
{\it generalized spectral shift function corresponding to the pair} $(A,B)$. 

A similar problem for unitary operators was considered by Neidhardt in \cite{N}. Let $U$ and $V$ be
unitary operators such that $V-U\in\bS_2$. Then $V=\exp({\rm i}A)U$, where $A$ is a self-adjoint operator
in $\bS_2$. Put $U_s=e^{{\rm i}sA}U$, $s\in\R$. It was shown in \cite{N} that there exists a function
$\eta\in L^1(\T)$ such that
\bay
\label{ntf}
\trace\left(\f(V)-\f(U)-\frac{d}{ds}\Big(\f(U_s)\Big)\Big|_{s=0}\right)=
\int_\T\f''\eta\,d\m,
\ey
whenever $\f''$ has absolutely convergent Fourier series. Such a function $\eta$ is unique modulo a
constant and it is called a {\it generalized spectral shift function corresponding to the pair} $(U,V)$.

We refer the reader to \cite{Ko2} and \cite{Bo} for applications of Koplienko's trace formula
\cite{Ko}.

In this paper we obtain better sufficient conditions on functions $\f$, under which trace formulae
\rf{ktf} and \rf{ntf} hold. We consider the case of unitary operators in \S 3 and the case of
self-adjoint operators in \S 4. We show that formula \rf{ktf} holds under the assumption
that $\f$ belongs to the Besov class $B_{\be1}^2(\R)$ while trace formula \rf{ntf} holds whenever
$\f\in B_{\be1}^2$. Note however, that the case of self-adjoint operators is considerably more
complicated. First of all, unlike in the case of functions on $\T$ the set of rational functions with
poles off $\R$ is not dense in $B_{\be1}^2(\R)$. Second, functions in $\f\in B_{\be1}^2(\R)$ do not have
to be Lipschitz and it is not clear how to interpret each of the operators 
$$
\f(B)-\f(A)\quad\mbox{and}\quad\frac{d}{ds}\Big(\f(A+sK)\Big)\Big|_{s=0}.
$$
However, it is still possible to define their difference and show that the difference belongs to $\bS_1$.

In \S 2 we outline the theory of double operator integrals developed by Birman and Solomyak in 
\cite{BS1}, \cite{BS2}, and \cite{BS4}, and we define Besov classes and discuss their properties. 

\

\section{\bf Preliminaries}
\setcounter{equation}{0}

\

In this section we collect necessary information on double operator integrals and Besov classes.

\medskip

{\bf Double operator integrals.} The technique of double operator integrals developed by Birman and
Solomyak in \cite{BS1}, \cite{BS2}, and \cite{BS4} plays an important role in perturbation theory.
We give here a brief introduction in this theory and state several results that will be used in the
main part of this paper.

Let $(\X,E)$ and $(\Y,F)$ be spaces with spectral measures $E$ and $F$ on a Hilbert space $\h$. Double
operator integrals are objects of the form
\bay
\label{doi}
\int\limits_\X\int\limits_\Y\psi(x,y)\,d E(x)T\,dF(y),
\ey
where $T$ is an operator on $\h$. Certainly, one has to specify how to understand the expression in 
\rf{doi}. Let us first define double operator integrals for bounded functions $\psi$ and operators $T$
of Hilbert Schmidt class $\bS_2$. Consider the spectral measure $\E$ whose values are orthogonal
projections on the Hilbert space $\bS_2$, which is defined by
$$
\E(\D\times\L)T=E(\D)TF(\L),\quad T\in\bS_2,
$$ 
for $\D$ and $\L$ being measurable subsets of $\X$ and $\Y$. Then $\E$ extends to a spectral measure on 
$\X\times\Y$ and if $\psi$ is a bounded measurable function on $\X\times\Y$, by definition
$$
\int\limits_\X\int\limits_\Y\psi(x,y)\,d E(x)T\,dF(y)=
\left(\,\,\int\limits_{\X\times\Y}\psi\,d\E\right)T.
$$
Clearly,
$$
\left\|\int\limits_\X\int\limits_\Y\psi(x,y)\,dE(x)T\,dF(y)\right\|_{\bS_2}
\le\|\psi\|_{L^\be}\|T\|_{\bS_2}.
$$
If 
$$
\int\limits_\X\int\limits_\Y\psi(x,y)\,d E(x)T\,dF(y)\in\bS_1
$$
for every $T\in\bS_1$, we say that $\psi$ is a {\it Schur multiplier of} $\bS_1$. In this case
by duality the map
$$
T\mapsto\int\limits_\X\int\limits_\Y\psi(x,y)\,d E(x)T\,dF(y)
$$
extends to a bounded transformer on the space of bounded linear operators on $\h$.

Suppose now that $A$ is a self-adjoint operator on a Hilbert space $\h$ and $B=A+K$, where 
$K$ is a self-adjoint operator of class $\bS_2$, and let $\f$ be a Lipschitz function on $\R$,
then $\f(B)-\f(A)\in\bS_2$ and
\bay
\label{s}
\f(B)-\f(A)=\iint\limits_{\R\times\R}\frac{\f(x)-\f(y)}{x-y}\,dE_B(x)K\,dE_A(y),
\ey
where $E_A$ and $E_B$ are spectral measure of $A$ and $B$. Here we can define the
function $(\f(x)-\f(y))(x-y)^{-1}$ on the diagonal $\{(x,x):~x\in\R\}$ in an arbitrary way.

A similar formula holds for unitary operators $U$ and $V$ with $V-U\in\bS_2$:
\bay
\label{ru}
\f(V)-\f(U)=\iint\limits_{\T\times\T}\frac{\f(\z)-\f(\t)}{\z-\t}\,dE_V(\z)(V-U)\,dE_U(\t),
\ey
where $\f$ is a Lipschitz function on $\T$. Again, the right-hand side
of this formula does not depend on the values of the function $(\f(\z)-\f(\t))(\z-\t)^{-1}$ on the
diagonal \lb$\{(\z,\z):~\z\in\T\}$. We refer the reader to 
\cite{BS1}, \cite{BS2}, and \cite{BS4} for more detailed information on double operator integrals. We
also mention recent survey article \cite{BS5}.

It follows from the results of \cite{F} and \cite{Pe1}, \cite{Pe2} mentioned in the introduction, that
the conditions $\f\in C^1$ and $\f'\in L^\be$ do not  imply that the above double operator integrals
determine bounded linear operators on $\bS_1$. On the other hand, it follows from the results of
\cite{Pe1} and \cite{Pe2}, that for functions $\f$ in the Besov class $B_{\be1}^1$ on the circle
and for functions $\f$ in the Besov class $B_{\be1}^1(\R)$ on $\R$ the following estimates hold:
\bay
\label{lu}
\left\|\,\,\iint\limits_{\T\times\T}\frac{\f(\z)-\f(\t)}{\z-\t}\,dE_V(\z)(V-U)\,dE_U(\t)\right\|_{\bS_1}
\le\const\|\f\|_{B_{\be1}^1}\|V-U\|_{\bS_1}
\ey
and
$$
\left\|\,\,\iint\limits_{\R\times\R}\frac{\f(x)-\f(y)}{x-y}\,dE_B(x)K\,dE_A(y)\right\|_{\bS_1}
\le\const\|\f\|_{B_{\be1}^1(\R)}\|K\|_{\bS_1}.
$$

In their papers \cite{BS1}, \cite{BS2}, and \cite{BS4} Birman and Solomyak studied the problem of the
differentiability of the map $t\mapsto\f(A+sK)$ in the operator norm and obtained sufficient
conditions (a similar problems was also studied there in the case of functions of unitary operators).
Later their results were improved in \cite{Pe1} and \cite{Pe2}. 

In this paper we need only differentiability results in the norm of $\bS_2$.
Let $\f$ be a function in $C^1(\R)$ such that $\f'\in L^\be$. Suppose that $A$ is a self-adjoint
operator (not necessarily bounded) and $K$ is a self-adjoint operator of class $\bS_2$. Then
\bay
\label{ds}
\frac{d}{ds}\Big(\f(A+sK)\Big)\Big|_{s=0}=\iint\limits_{\R\times\R}\frac{\f(x)-\f(y)}{x-y}\,dE_A(x)K\,dE_A(y)
\ey
(the derivative exists in the $\bS_2$ norm). This follows from formula \rf{s} and Proposition 3.2 of
\cite{dPS}.

A similar result holds for functions of unitary operators. Let $\f\in C^1(\T)$. Suppose that $U$
is a unitary operator, $A$ is a self-adjoint operator of class $\bS_2$. Then
\bay
\label{du}
\frac{d}{ds}\Big(\f\big(e^{{\rm i}sA}U\big)\Big)\Big|_{s=0}=
{\rm i}\iint\limits_{\T\times\T}\t\frac{\f(\z)-\f(\t)}{\z-\t}\,dE_U(\z)A\,dE_U(\t).
\ey
The proof of this formula is much simpler than in the case of possibly unbounded self-adjoint operators.

\medskip

{\bf Besov classes.} Let $0<p,\,q\le\be$ and $s\in\R$. The Besov class $B^s_{pq}$ of functions (or
distributions) on $\T$ can be defined in the following way. Let $v$ be a $C^\be$ function on $\R$ such
that
\bay
\label{v}
w\ge0,\quad\supp w\subset\left[\frac12,2\right],\quad\mbox{and}\quad
\sum_{n=-\be}^\be w(2^{n}x)=1\quad\mbox{for}\quad x>0.
\ey
Consider the trigonometric polynomials $W_n$, and $W_n^\#$ defined by
$$
W_n(z)=\sum_{n\in\Z}v\left(\frac{k}{2^n}\right),\quad n\ge1,\quad W_0(z)=\bar z+1+z,\quad
\mbox{and}\quad W_n^\#(z)=\ov{W_n(z)},\quad n\ge0.
$$
Then for each distribution $\f$ on $\T$
$$
\f=\sum_{n\ge0}\f*W_n+\sum_{n\ge1}\f*W^\#_n.
$$
The Besov class $B^s_{pq}$ consists of functions (in the case $s>0$) or distributions $\f$ on $\T$
such that
$$
\{\|2^{ns}\f*W_n\|_{L^\p}\}_{n\ge0}\in\ell^q\quad\mbox{and}
\quad\{\|2^{ns}\f*W^\#_n\|_{L^\p}\}_{n\ge0}\in\ell^q
$$
Besov classes admit many other descriptions. in particular, for $s>0$ the space $B^s_{pq}$ can be
described in terms of moduli of continuity (or moduli of smoothness).

To define (homogeneous) Besov classes $B^s_{pq}(\R)$ on the real line, we consider the same function $w$
as in \rf{v} and define the functions $W_n$ and $W^\#_n$ on $\R$ by
$$
\F W_n(x)=w\left(\frac{x}{2^n}\right),\quad\F W^\#_n(x)=\F W_n(-x),\quad n\in\Z,
$$
where $\F$ is the {\it Fourier transform}. The Besov class $B^s_{pq}(\R)$ consists of
distributions $\f$ on $\R$ such that
$$
\{\|2^{ns}\f*W_n\|_{L^\p}\}_{n\in\Z}\in\ell^q(\Z)\quad\mbox{and}
\quad\{\|2^{ns}\f*W^\#_n\|_{L^\p}\}_{n\in\Z}\in\ell^q(\Z)
$$
According to this definition, the space $B^s_{pq}(\R)$ contains all polynomials. However, it is not
necessary to include all polynomials. 

In this paper we need only Besov spaces $B_{\be1}^1$ and $B_{\be1}^2$. In the case of functions
on the real line it is convenient to restrict the degree of polynomials in $B_{\be1}^1(\R)$ by $1$
and in $B_{\be1}^2(\R)$ by $2$. It is also convenient to consider the following seminorms on
$B_{\be1}^1(\R)$ and in $B_{\be1}^2(\R)$:
$$
\|\f\|_{B_{\be1}^1(\R)}=\sup_{x\in\R}|\f'(x)|+\sum_{n\in\Z}2^n\|\f*W_n\|_{L^\be}+
\sum_{n\in\Z}2^n\|\f*W^\#_n\|_{L^\be}
$$
and
$$
\|\f\|_{B_{\be1}^2(\R)}=\sup_{x\in\R}|\f''(x)|+\sum_{n\in\Z}2^{2n}\|\f*W_n\|_{L^\be}+
\sum_{n\in\Z}2^{2n}\|\f*W^\#_n\|_{L^\be}.
$$
The classes $B_{\be1}^1(\R)$ and $B_{\be1}^2(\R)$ can be described as classes of function on $\R$ 
in the following way:
$$
\f\in B_{\be1}^1(\R)\quad\Longleftrightarrow\quad
\sup_{t\in\R}|\f'(t)|+\int\limits_\R\frac{\|\D_t^2\f\|_{L^\be}}{|t|^2}\,dt<\be
$$
and
$$
\f\in B_{\be1}^2(\R)\quad\Longleftrightarrow\quad
\sup_{t\in\R}|\f''(t)|+\int\limits_\R\frac{\|\D_t^3\f\|_{L^\be}}{|t|^3}\,dt<\be,
$$
where $\D_t$ is the difference operator defined by $(\D_t\f)(x)=\f(x+t)-\f(x)$.

It is interesting to note that the Besov class $B_{\be1}^2(\R)$ also appears in a natural way
in perturbation theory in \cite{Pe3}, where the following problem is studied: in which case
$$
\f(T_f)-T_{\f\circ f}\in \bS_1?
$$ 
($T_g$ is a Toeplitz operator with symbol $g$.)

We refer the reader to \cite{Pe} for more detailed information on Besov classes.

\

\section{\bf The case of unitary operators}
\setcounter{equation}{0}

\

Let $U$ and $V$ be unitary operators such that $V-U\in \bS_2$. Denote by $E_U$ and $E_V$
the spectral measures of $U$ and $V$. Let $A$ be a self-adjoint operator such that 
$\s(A)\subset[-\pi,\pi]$ and $V=\exp({\rm i}A)U$. It is easy to see that $A\in\bS_2$.

Put
\bay
\label{us}
U_s=e^{{\rm i}sA}U.
\ey

Consider the class $\Li\,\hat\odot\,L^\be$ that consists of functions $u$ on $\T\times\T$
that admit a representation
\bay
\label{tens}
u(\z,\t)=\sum_{n\ge0}f_n(\z)g_n(\t),\quad\z,\,\t\in\T,
\ey
where $f_n\in\Li$, $g_n\in L^\be$, and
\bay
\label{proj}
\sum_{n\ge0}\|f_n\|_{\Li}\cdot\|g_n\|_\be<\be.
\ey
By definition, $\|u\|_\te$ is the infimum of the left-hand side of \rf{proj} over all
functions $f_n$ and $g_n$ satisfying \rf{tens}. We consider here the following seminorm on the space
$\Li$ of Lipschitz functions:
$$
\|f\|_\Li=\sup_{\z\ne\t}\frac{|f(\z)-f(\t)|}{|\z-\t|}.
$$

For a differentiable function $\f$ on $\T$ we define the function $\breve\f$ 
on $\T\times\T$ by
$$
\breve\f(\z,\t)=\left\{\begin{array}{ll}
\frac{\f(\z)-\f(\t)}{\z-\t},&\z\ne\t,\\[.3cm]
\f'(\z),&\z=\t.\end{array}\right.
$$

\begin{thm}
\label{S2}
If $\f\in B_{\infty1}^2$, then
\bay
\label{s1}
\f(V)-\f(U)-\frac{d}{ds}\Big(\f(U_s)\Big)\Big|_{s=0}\in\bS_1
\ey
and
$$
\left\|\f(V)-\f(U)-\frac{d}{ds}\Big(\f(U_s)\Big)\Big|_{s=0}\right\|_{\bS_1}
\le\const\|\f\|_{B_{\infty1}^2}\|V-U\|_{\bS_2}.
$$
\end{thm}

To prove Theorem \ref{S2}, we need the following fact.

\begin{thm}
\label{bes}
If $\f\in B_{\infty1}^2$, then $\breve\f\in\te$
and
$$
\|\breve\f\|_\te\le\const\|\f\|_{B_{\infty1}^2}.
$$
\end{thm}

\Pf We have 
\bay
\label{+-}
\breve\f(\z,\t)=\sum_{j,k\ge0}\hat\f(j+k+1)\z^j\t^k+
\sum_{j,k<0}\hat\f(j+k+1)\z^j\t^k.
\ey

Let us show that the first term on the right-hand side of \rf{+-} belongs to $\te$.
A similar result for the second term in \rf{+-} can
be proved in the same way. We use the construction given in the proof of Theorem 2 of Section 3 of
\cite{Pe1}. We have
\bay
\label{ab}
\sum_{j,k\ge0}\!\hat\f(j+k+1)\z^j\t^k\!=\!\sum_{j,k\ge0}\!\a_{jk}\hat\f(j+k+1)\z^j\t^k
\!+\!\sum_{j,k\ge0}\!\b_{jk}\hat\f(j+k+1)\z^j\t^k\!,
\ey
where 
$$
\a_{jk}=\left\{\begin{array}{ll}\frac12,&j=k=0,\\[.2cm]
\frac{2j-k}{j+k},&j+k>0,~\frac k2\le j\le2k,\\[.2cm]
0,&j\ge2k,
\end{array}\right.\qquad\mbox{and}\qquad\b_{jk}=1-\a_{jk}.
$$

Let us prove that the function
$$
(\z,\t)\mapsto\sum_{j,k\ge0}\!\b_{jk}\hat\f(j+k+1)\z^j\t^k
$$
on the right-hand side of \rf{ab} belongs to $\te$. 

We define the functions $q$ and $r$ on $\R$ by
\bay
\label{rq}
q(x)=\left\{\begin{array}{ll}0,&x\le\frac12,\\[.2cm]
\frac{2x-1}{x+1},& \frac12\le x\le2,\\[.2cm]
1,&x\ge2,
\end{array}\right.\qquad\mbox{and}\qquad
r(x)=\left\{\begin{array}{ll}0,&x\le\frac32,\\[.2cm]
\frac{2x-3}{x},& \frac32\le x\le3,\\[.2cm]
1,&x\ge3.
\end{array}\right.
\ey
Put
$$
Q_n(z)=\sum_{j\ge0}q\left(\frac jn\right),\quad R_n(z)=\sum_{j\ge0}r\left(\frac jn\right)
\quad\mbox{for}\quad n>0
$$
and
$$
Q_0(z)=R_0(z)=\frac12+\sum_{j\ge1}z^j.
$$
It is easy to see that
\bay
\label{b}
\sum_{j,k\ge0}\!\b_{jk}\hat\f(j+k+1)\z^j\t^k=
\sum_{n\ge0}\z^n\Big(\big((S^*)^n\psi*Q_n\big)(\t)\Big),
\ey
where $\psi=\pp_+\bar z\f$ and $S^*\psi=\frac{\psi-\psi(0)}{z}$. We have
$$
\sum_{n\ge0}\|z^n\|_{\Li}\|(S^*)^n\psi*Q_n\|_\be\le\const\sum_{n\ge0}n\|(S^*)^n\psi*Q_n\|_\be
=\const\sum_{n\ge0}n\|\psi*R_n\|_\be.
$$

Let us show that for $\f\in B_{\infty1}^2$,
$$
\sum_{n\ge0}n\|\psi*R_n\|_\be<\be.
$$

Consider the function $r^\flat$ on $\R$ defined by $r^\flat(x)=1-r(|x|)$, $x\in\R$. Put
$$
R^\flat_n(z)=\sum_{j\in\Z}r^\flat\left(\frac jn\right),\quad n>0.
$$
then $\|R^\flat_n\|_{L^1}\le\const$ (see the proof of Lemma 3 of \cite{Pe1}). Suppose that
$n\ge 2^{m}$. Then
$$
R_n*\psi=R_n*\sum_{k\ge m}\psi*W_k
$$
(see \S 2 for the definition of the polynomials $W_k$). Hence, 
$$
\|R_n*\psi\|_\be\le\left\|R_n*\left(\sum_{k\ge m}\psi*W_k\right)\right\|_\be
\le\big(1+\|R^\flat_n\|_1\big)\sum_{k\ge m}\left\|\psi*W_k\right\|_\be.
$$
It follows that 
\begin{align*}
\sum_{n\ge2}n\|\psi*R_n\|_\be&=\sum_{m\ge1}~\sum_{n=2^m}^{2^{m+1}-1}n\|\psi*R_n\|_\be\\[.2cm]
&\le\const\sum_{m\ge1}2^{2m}\sum_{k\ge m}\|\psi*W_k\|_\be
\le\const\sum_{m\ge1}2^{2m}\|\psi*W_m\|_\be<\be,
\end{align*}
since $\psi\in B_{\be1}^2$.

Let us now show that the function
$$
(\z,\t)\mapsto\sum_{j,k\ge0}\!\a_{jk}\hat\f(j+k+1)\z^j\t^k
$$
belongs to the space $\te$.

It follows from \rf{b} that 
$$
\sum_{j,k\ge0}\!\a_{jk}\hat\f(j+k+1)\z^j\t^k=
\sum_{n\ge0}\Big(\big((S^*)^n\psi*Q_n\big)(\z)\Big)\t^n.
$$
It suffices to show that
$$
\sum_{n\ge0}\big\|(S^*)^n\psi*Q_n\big\|_{\Li}<\be.
$$
By the Bernstein inequality, we have
\ba
\sum_{n\ge0}\big\|(S^*)^n\psi*Q_n\big\|_{\Li}&=
\sum_{n\ge0}\Big\|\Big((S^*)^n\psi*Q_n\Big)'\Big\|_\be\\[.2cm]
&\le\sum_{n\ge0}\sum_{k\ge0}\Big\|\Big(\big((S^*)^n(\psi*W_k)\big)*Q_n\Big)'\Big\|_\be\\[.2cm]
&\le\sum_{n\ge0}\sum_{k\ge0}2^{k+1}\big\|\big((S^*)^n(\psi*W_k)\big)*Q_n\big\|_\be\\[.2cm]
&\le\sum_{n\ge0}\sum_{k\ge0}2^{k+1}\big\|\psi*W_k*R_n\big\|_\be\\[.2cm]
&\le\sum_{0\le n\le2^{k+2}/3}\big(1+\|R^\flat_n\|_1\big)\sum_{k\ge0}2^{k+1}\|\psi*W_k\|_\be\\[.2cm]
&\le\const\sum_{k\ge0}2^{2k}\|\psi*W_k\|_\be\le\const\|\psi\|_{B_{\be1}^2},
\end{align*}
since, clearly, $\psi*W_k*R_n=\0$ if $n>2^{k+2}/3$. $\bl$


\medskip

{\bf Proof of Theorem \ref{S2}.} Since $\f\in C^1(\T)$, we have by \rf{du}
$$
\frac{d}{ds}\Big(\f(U_s)\Big)\Big|_{s=0}=
{\rm i}\iint\limits_{\T\times\T}\t\breve\f(\z,\t)\,dE_U(\z)A\,dE_U(\t).
$$
On the other hand, by \rf{ru},
\begin{align*}
\f(V)-\f(U)&=\iint\limits_{\T\times\T}\breve\f(\z,\t)\,dE_V(\z)(V-U)\,dE_U(\t)\\[.2cm]
&=-\iint\limits_{\T\times\T}\t\breve\f(\z,\t)\,dE_V(\z)(I-VU^*)\,dE_U(\t)\\[.2cm]
&=-\iint\limits_{\T\times\T}\t\breve\f(\z,\t)\,dE_V(\z)
\big(I-e^{{\rm i}A}\big)\,dE_U(\t).
\end{align*}
Thus
\begin{align*}
\f(V)-\f(U)-\frac{d}{ds}\Big(\f(U_s)\Big)\Big|_{s=0}
&=
-\iint\limits_{\T\times\T}\t\breve\f(\z,\t)\,dE_V(\z)
\big(I-e^{{\rm i}A}\big)\,dE_U(\t)\\[.2cm]
&-{\rm i}\iint\limits_{\T\times\T}\t\breve\f(\z,\t)\,dE_U(\z)A\,dE_U(\t)\\[.2cm]
&=-\iint\limits_{\T\times\T}\t\breve\f(\z,\t)\,dE_V(\z)
\big(I-e^{{\rm i}A}\big)\,dE_U(\t)\\[.2cm]
&+\iint\limits_{\T\times\T}\t\breve\f(\z,\t)\,dE_U(\z)
\big(I-e^{{\rm i}A}\big)\,dE_U(\t)\\[.2cm]
&+\iint\limits_{\T\times\T}\t\breve\f(\z,\t)\,dE_U(\z)
\big(e^{{\rm i}A}-I-{\rm i}A\big)\,dE_U(\t).
\end{align*}

It is easy to see that $e^{{\rm i}A}-I-{\rm i}A\in\bS_1$, and so by \rf{lu},
$$
\iint\limits_{\T\times\T}\t\breve\f(\z,\t)\,dE_U(\z)
\big(e^{{\rm i}A}-I-{\rm i}A\big)\,dE_U(\t)\in \bS_1
$$
and
$$
\left\|\iint\limits_{\T\times\T}\t\breve\f(\z,\t)\,dE_U(\z)
\big(e^{{\rm i}A}-I-{\rm i}A\big)\,dE_U(\t)\right\|_{\bS_1}\le
\const\|\f\|_{B_{\infty1}^1}\|\big(e^{{\rm i}A}-I-{\rm i}A\big)\|_{\bS_1}.
$$
Clearly, $\|\f\|_{B_{\infty1}^1}\le\const\|\f\|_{B_{\infty1}^2}$.

On the other hand, let $\{f_n\}_{n\ge0}$ and
$\{g_n\}_{n\ge0}$ be sequences of functions such that
$$
\breve\f(\z,\t)=\sum_{n\ge0}f_n(\z)g_n(\t),\quad\z,\,\t\in\T,
$$
and \rf{proj} holds.
We have 
\begin{align*}
&\iint\limits_{\T\times\T}\t\breve\f(\z,\t)\,dE_V(\z)
\big(I-e^{{\rm i}A}\big)\,dE_U(\t)-
\iint\limits_{\T\times\T}\t\breve\f(\z,\t)\,dE_U(\z)
\big(I-e^{{\rm i}A}\big)\,dE_U(\t)\\[.2cm]
=&\iint\limits_{\T\times\T}\sum_{n\ge0}f_n(\z)\t g_n(\t)\,dE_V(\z)
\big(I-e^{{\rm i}A}\big)\,dE_U(\t)\\[.2cm]
&-\iint\limits_{\T\times\T}\sum_{n\ge0}f_n(\z)\t g_n(\t)\,dE_U(\z)
\big(I-e^{{\rm i}A}\big)\,dE_U(\t)\\[.2cm]
=&\sum_{n\ge0}f_n(V)\big(I-e^{{\rm i}A}\big)g_n(U)U-
\sum_{n\ge0}f_n(U)\big(I-e^{{\rm i}A}\big)g_n(U)U\\[.2cm]
=&\sum_{n\ge0}\big(f_n(V)-f_n(U)\big)\big(I-e^{{\rm i}A}\big)g_n(U)U.
\end{align*}
Thus
\begin{align*}
&\left\|\,\iint\limits_{\T\times\T}\t\breve\f(\z,\t)\,dE_V(\z)
\big(I-e^{{\rm i}A}\big)\,dE_U(\t)-
\iint\limits_{\T\times\T}\t\breve\f(\z,\t)\,dE_U(\z)
\big(I-e^{{\rm i}A}\big)\,dE_U(\t)\right\|_{\bS_1}\\[.2cm]
&\le\sum_{n\ge0}\big\|\big(f_n(V)-f_n(U)\big)\big\|_{\bS_2}
\big\|\big(I-e^{{\rm i}A}\big)\big\|_{\bS_2}\|g_n(U)\|\\[.2cm]
&\le\const\big\|\big(I-e^{{\rm
i}A}\big)\big\|_{\bS_2}\sum_{n\ge0}\|f_n\|_\Li\cdot\|g_n\|_{\be}<\be.
\end{align*}
This completes the proof. $\bl$

Let now $\eta$ be a generalized spectral shift function for the pair $(V,U)$.

\begin{thm}
\label{trf}
Let $U$ and $V$ be unitary operators such that $V-U\in\bS_2$ and let $U_s$ be defined by 
{\em\rf{trf}}. Then for any $\f\in B^2_{\be1}$,
\bay
\label{tfu}
\trace\left(\f(V)-\f(U)-\frac{d}{ds}\Big(\f(U_s)\Big)\Big|_{s=0}\right)=
\int_\T\f''\eta\,d\m.
\ey
\end{thm}

\medskip

\Pf Since clearly, $B_{\be1}^2\subset C^1$, the fact that the operator in \rf{s1}
belongs to $\bS_1$ is an immediate consequence of Theorems \ref{S2} and \ref{bes}.

It is easy to see from the definition of the space $B^2_{\be1}$ given in \S 2 that the
trigonometric polynomials are dense in $B^2_{\be1}$. Let $\f_n$ be trigonometric
polynomials such that 
$$
\lim_{n\to\be}\|\f-\f_n\|_{B^2_{\be1}}=0.
$$
Since $B^2_{\be1}$ is continuously imbedded in the space $C^2$ of functions with two
continuous derivatives, it follows that $\f_n\to\f$ in $C^2$. Since $\eta\in L^1$, it follows
that
$$
\lim_{n\to\be}\int_\T\f_n''\eta\,d\m=\int_\T\f''\eta\,d\m.
$$
On the other hand, it follows from Theorems \ref{S2} and \ref{bes} that 
$$
\left\|\left(\f_n(V)-\f_n(U)-
\frac{d}{ds}\Big(\f_n(U_s)\Big)\Big|_{s=0}\right)
-\left(\f(V)-\f(U)-\frac{d}{ds}\Big(\f(U_s)\Big)\Big|_{s=0}\right)\right\|_{\bS_1}\to0
$$
as $n\to\be$.
The result follows now from the fact that trace formula \rf{tfu}
is valid for all trigonometric polynomials $\f$ (see \cite{N}). $\bl$

\

\section{\bf The case of self-adjoint operators}
\setcounter{equation}{0}

\

In this section we extend Koplienko's trace formula for self-adjoint operators to a considerably bigger
class of functions. 

Let $A$ be a self-adjoint operator (not necessarily bounded) on Hilbert space and let $K$ be a
self-adjoint operator of class $\bS_2$. Put $B=A+K$. As we have already mentioned in the
introduction, Koplienko introduced in \cite{Ko} the generalized spectral shift function $\eta\in L^1$ that
corresponds to the pair $(A,B)$ and showed that for rational functions $\f$ with poles off the real line
the following trace formula holds. 

\bay
\label{trf}
\trace\left(\f(B)-\f(A)-\frac{d}{ds}\Big(\f(A_s)\Big)\Big|_{s=0}\right)=\int_\R\f''(x)\eta(x)\,dx.
\ey

We are going to extend this formula to the Besov class $B_{\be1}^2(\R)$. Note however, that 
the situation with self-adjoint operators is subtler than with unitary operators. First of all, the
rational functions are not dense in $B_{\be1}^2(\R)$ and this makes it more difficult to extend
formula \rf{trf} from rational functions to $B_{\be1}^2(\R)$. Secondly,
functions  in $B_{\be1}^2(\R)$ do not have to belong to the space $\Li$ of Lipschitz functions
on $\R$, which we equip with the seminorm:
$$
\|f\|_\Li=\sup_{x\ne y}\frac{|f(x)-f(y)|}{|x-y|}.
$$
Thus for $\f\in B_{\be1}^2(\R)$, none of the operators
$$
\f(B)-\f(A)\quad\mbox{and}\quad\frac{d}{ds}\Big(\f(A_s)\Big)\Big|_{s=0}
$$
has to be in $\bS_2$.  In fact, it is not clear how one can interpret each of those operators.
However, it turns out that their difference still makes sense for functions $\f\in B_{\be1}^2(\R)$
and formula \rf{trf} holds for such functions $\f$.

To do it, we first prove formula \rf{trf} in the case $\f\in B_{\be1}^2(\R)\bigcap\Li$ and estimate
the $\bS_1$ norm of the left-hand side of \rf{trf} in terms of $\|\f\|_{B_{\be1}^2}$. Then we define
the operator on the left-hand side of \rf{trf} for functions $f\in B_{\be1}^2(\R)$ and prove formula
\rf{trf} for such functions.

For a differentiable function $\f$ on $\R$ we define the function $\breve\f$ on
$\R\times\R$ by
$$
\breve\f(x,y)=\left\{\begin{array}{ll}
\frac{\f(x)-\f(y)}{\z-\t},&x\ne y,\\[.3cm]
\f'(x),&x=y.\end{array}\right.
$$

We consider in this section the space $\ite$  of functions $u$ on $\R\times\R$ that admit a
representation
\bay
\label{ir}
u(x,y)=\int_\O f_\o(x)g_\o(y)\,d\mu(\o),
\ey
where $(\O,\mu)$ is a measure space and the functions $(\o,x)\mapsto f_\o(x)$ and
$(\o,y)\mapsto g_\o(y)$ are measurable functions on $\O\times\R$ such that $f_\o\in\Li$,
$g_\o\in L^\be$ for almost all $\o\in\O$, and
\bay
\label{tn}
\int_\O\|f_\o\|_\Li\cdot\|g_\o\|_{L^\be}\,d\mu(\o)<\be.
\ey
By definition, the norm of $u$ in $\ite$ is the infimum of the left-hand side of \rf{tn} over
all representations of the form \rf{ir}.

\begin{thm}
\label{M}
Let $M>0$. Suppose that $\f$ is a bounded function on $\R$ such that $\supp\F\f\subset[M/2,2M]$.
Then 
\bay
\label{pdp}
\f(B)-\f(A)-\frac{d}{ds}\Big(\f(A_s)\Big)\Big|_{s=0}\in\bS_1
\ey
and
\bay
\label{MK}
\left\|\f(B)-\f(A)-\Big(\frac{d}{ds}\f(A_s)\Big)\Big|_{s=0}\right\|_{\bS_1}
\le\const\cdot M^2\|K\|_{\bS_2}^2\|\f\|_{L^\be}.
\ey

\end{thm}

To prove Theorem \ref{M}, we need the following fact.

\begin{lem}
\label{L1}
Let $\f$ be a function on $\R$ such that $\supp\F\f\subset[M/2,2M]$.
Then \lb$\breve\f\in\ite$ and 
$$
\|\breve\f\|_\ite\le\const\cdot M^2\|\f\|_{L^\be}.
$$
\end{lem}

\Pf Let $q$ and $r$ be the functions on $\R$ defined by \rf{rq}.

Consider the distributions $Q_t$ and $R_t$, $t>0$, on $\R$ such that
$$
(\F Q_t)(x)=q(x/t)\quad\mbox{and}\quad(\F R_t)(x)=r(x/t).
$$

It was shown in \cite{Pe2} (formula (5)) that
\bay
\label{gf}
\breve\f(x,y)=\int_0^\be\Big(\big(S^*_t\f\big)*Q_t\Big)(x)e^{{\rm i}ty}\,dt+
\int_0^\be\Big(\big(S^*_t\f\big)*Q_t\Big)(y)e^{{\rm i}tx}\,dt,
\ey
where $(S_t^*\f)(x)=e^{-{\rm i}tx}\f(x)$.

Clearly, 
$$
\|\breve\f\|_\te\le\int_0^\be\big\|\big(S^*_t\f\big)*Q_t\big\|_\Li\,dt
+\int_0^\be\big\|\big(S^*_t\f\big)*Q_t\big\|_{L^\be}t\,dt.
$$

By the Bernstein inequality,
$$
\big\|(S^*_t\f\big)*Q_t\big\|_\Li=\Big\|\Big(\big(S^*_t\f\big)*Q_t\Big)'\Big\|_{L^\be}
\le2M\big\|\big(S^*_t\f\big)*Q_t\big\|_{L^\be}, 
$$
and so
\ba
\int_0^\be\big\|\big(S^*_t\f\big)*Q_t\big\|_\Li\,dt&\le
2M\int_0^\be\big\|\big(S^*_t\f\big)*Q_t\big\|_{L^\be}\,dt
=2M\int_0^\be\big\|\f*R_t\big\|_{L^\be}\,dt\\[.2cm]
&=2M\int_0^{4M/3}\big\|\f*R_t\big\|_{L^\be}\,dt
\le\frac83M^2\big\|\f*R_t\big\|_{L^\be},
\end{align*}
since, obviously, $(S^*_t\f)*R_t=\0$ for $t\ge4M/3$.

On the other hand, 
\ba
\int_0^\be\big\|\big(S^*_t\f\big)*Q_t\big\|_{L^\be}t\,dt&
=\int_0^\be\big\|\f*R_t\big\|_{L^\be}t\,dt\\[.2cm]
&=\int_0^{4M/3}\big\|\f*R_t\big\|_{L^\be}t\,dt.
\end{align*}

It remains to observe that if $R_t^\flat$ is the function on $\R$ such that
$$
\F\big(R_t^\flat\big)(x)=1-\F(R_t)(|x|),
$$
then $R_t^\flat\in L^1$, $\|R_t^\flat\|_{L^1}$ does not depend on $t$ and
$$
\big\|\f*R_t\big\|_{L^\be}\le(1+\big\|R_t^\flat\big\|_{L^1})\|\f\|_{L^\be}.\quad\bl
$$

\medskip

{\bf Proof of Theorem \ref{M}.} By \rf{s} and \rf{ds}, we have 
\begin{align}
\label{doif}
\f(B)-\f(A)-\frac{d}{ds}\Big(\f(A_s)\Big)\Big|_{s=0}&=
\iint\limits_{\R\times\R}\breve\f(x,y)dE_B(x)KdE_A(y)\nonumber\\[.2cm]
&-\iint\limits_{\R\times\R}\breve\f(x,y)dE_A(x)KdE_A(y).
\end{align} 

By Lemma \ref{L1}, $\breve\f$ admits a representation
$$
\breve\f(x,y)=\int_\O f_\o(x)g_\o(y)\,d\mu(\o)
$$
such that
$$
\int_\O\|f_\o\|_\Li\cdot\|g_\o\|_{L^\be}\,d\mu(\o)\le\const\cdot M^2\|\f\|_{L^\be}.
$$
We have
\ba
\iint\limits_{\R\times\R}\breve\f(x,y)dE_B(x)KdE_A(y)&=
\int\limits_\O\left(\,\,\iint\limits_{\R\times\R}f_\o(x)g_\o(y)dE_B(x)KdE_A(y)\right)d\mu(\o)
\\[.2cm]
&=\int\limits_\O f_\o(B)Kg_\o(A)\,d\mu(\o).
\end{align*}
Similarly,
$$
\iint\limits_{\R\times\R}\breve\f(x,y)\,dE_A(x)K\,dE_A(y)
=\int\limits_\O f_\o(A)Kg_\o(A)d\mu(\o).
$$
Thus
$$
\f(B)-\f(A)-\frac{d}{ds}\Big(\f(A+sK)\Big)\Big|_{s=0}=
\int\limits_\O \big(f_\o(B)-f_\o(A)\big)Kg_\o(A)\,d\mu(\o).
$$
Hence,
\ba
\left\|\f(B)-\f(A)-\frac{d}{ds}\Big(\f(A_s)\Big)\Big|_{s=0}\right\|_{\bS_1}&\le
\int\limits_\O\big\|f_\o(B)-f_\o(A)\big\|_{\bS_2}\|K\|_{\bS_2}\|g_\o(A)\|\,d\mu(\o)\\[.2cm]
&\le\|K\|_{\bS_2}\int\limits_\O\|f_\o\|_\Li\|B-A\|_{\bS_2}\|g_\o\|_{L^\be}\,d\mu(\o)\\[.2cm]
&=\|K\|^2_{\bS_2}\int\limits_\O\|f_\o\|_\Li\|g_\o\|_{L^\be}\,d\mu(\o)\\[.2cm]
&\le\const\cdot M^2\|K\|^2_{\bS_2}\|\f\|_{L^\be}.\quad\bl
\end{align*}

\medskip

{\bf Remark.} It is easy to see that the same conclusion holds if $\f$ is a bounded function on
$\R$ such that $\supp\F\f\subset[-2M,-M/2]$.

\begin{thm}
\label{int}
Suppose that $\f\in B^2_{\be1}(\R)\bigcap\Li$. Then {\em\rf{pdp}} holds,
$$
\left\|\f(B)-\f(A)-\frac{d}{ds}\Big(\f(A_s)\Big)\Big|_{s=0}\right\|_{\bS_1}
\le\const\cdot \|K\|_{\bS_2}^2\|\f\|_{B^2_{\be1}},
$$
and {\em\rf{trf}} holds.
\end{thm}

We need the following lemma.

\begin{lem}
\label{uL}
Let $\{f_n\}_{n\ge0}$ and $f$ be functions in $\Li(\R)$ such that
$$
\lim_{n\to\be}f_n(x)=f(x),\quad x\in\R,\quad\mbox{and}\quad \sup_n\|f_n\|_\Li<\be.
$$
Then
$$
\lim_{n\to\be}\big(f_n(B)-f_n(A)\big)=f(B)-f(A)
$$
in $\bS_2$.
\end{lem}

Let us first prove Theorem \ref{int}.

\medskip

{\bf Proof of Theorem \rf{int}.} 
Since $\f\in B^2_{\be1}(\R)$, $\f$ is continuously differentiable,
and so both operators 
$$
\f(B)-\f(A)\quad\mbox{and}\quad\frac{d}{ds}\Big(\f(A+sK)\Big)\Big|_{s=0}
$$
belong to $\bS_2$ (see \S 2).

Clearly, if $\f$ is a linear function, then the operator in \rf{pdp}
is zero. Suppose first that $\F\f''\in L^1$. Then
$$
\f=\sum_{n\in\Z}\big(\f_n+\f^\#_n\big).
$$
where 
$$
\f_n=\f*\F^{-1}\chi_{[2^n,2^{n+1}]}\quad\mbox{and}\quad\f^\#_n=\f*\F^{-1}\chi_{[-2^{n+1},-2^n]}.
$$
Clearly,
\bay
\label{fn}
2^{2n}\|\f_n\|_{L^\be}\le\const\|\F\f''_n\|_{L^1}\quad\mbox{and}\quad
2^{2n}\|\f^\#_n\|_{L^\be}\le\const\|\F\f^\#_n\|_{L^1}.
\ey
By Theorem \ref{M},
$$
\sum_{n\in\Z}\left\|\f_n(B)-\f_n(A)-\frac{d}{ds}\Big(\f_n(A_s)\Big)\Big|_{s=0}\right\|_{\bS_1}
\le\const\sum_{n\in\Z}2^{2n}\|\f_n\|_{L^\be}
$$
and the same estimate holds for the functions $\f_n^\#$ in place of $\f_n$.
It follows now from \rf{fn} that
\ba
\left\|\f(B)-\f(A)-\frac{d}{ds}\Big(\f(A_s)\Big)\Big|_{s=0}\right\|_{\bS_1}&\le
\const\sum_{n\in\Z}2^{2n}\Big(\|\f_n\|_{L^\be}+\|\f^\#_n\|_{L^\be}\Big)\\[.2cm]
&\le\const\|\F\f''\|_{L^1}.
\end{align*}
Since the rational functions are dense in the space $\{\f:~\F\f''\in L^1\}$ and trace formula
\rf{trf} holds for rational functions with poles outside $\R$ (Koplienko's theorem \cite{Ko}),
it is easy to see that it also holds for arbitrary functions $\f$ with $\F\f''\in L^1$.

Suppose now that $\f\in B^2_{\be1}$. Since 
$$
\sum_{n\in\Z}2^{2n}\Big(\|\f*W_n\|_{L^\be}+\|\f*W_n^\#\|_{L^\be}\Big)<\be,
$$
and inequality \rf{MK} holds,
it suffices to show that formula \rf{trf} holds for functions $\f*W_n$ and $\f*W_n^\#$.

The following argument is similar to the argument given in the proof of Theorem 4 of \cite{Pe2}
to establish the Lifshitz--Krein trace formula for functions in $B^1_{\be1}(\R)$.
Put $\psi=\f*V_n$. Then $\supp\psi\subset\big[2^{n-1},2^{n+1}\big]$. Consider a smooth
nonnegative function $h$ on $\R$ such that $\supp h\subset[-1,1]$ and $\int_{-1}^1h(x)\,dx=1$.
For $\e>0$ put $h_\e(x)=\e^{-1}h(x/\e)$. Let $\psi_\e$ be the function defined by
$\F\psi_\e=\F\psi*h_\e$. Clearly, 
$$
\F\psi_\e\in L^1,\quad\lim_{\e\to0}\|\psi_\e\|_{L^\be}=\|\psi\|_{L^\be},\quad\mbox{and}\quad
\lim_{\e\to0}\psi_\e(x)=\psi(x)\quad\mbox{for}\quad x\in\R.
$$
Then formula \rf{trf} holds for $\psi_\e$. Clearly, 
$$
\lim_{\e\to0}\int_\R\psi''_\e(x)\eta(x)\,dx=\int_\R\psi''(x)\eta(x)\,dx.
$$
Thus to prove that \rf{trf} holds for $\psi$, it suffices to show that
$$
\lim_{\e\to0}\trace\!\left(\!\psi_e(B)-\psi_\e(A)
-\frac{d}{ds}\Big(\psi_\e(A_s)\Big)\Big|_{s=0}\!\right)\!
=\trace\!\left(\!\psi(B)-\psi(A)-\frac{d}{ds}\Big(\psi(A_s)\Big)\Big|_{s=0}\!\right)\!.
$$
By \rf{doif}, we have 
\ba
\psi_e(B)-\psi_\e(A)-\frac{d}{ds}\Big(\psi_\e(A_s)\Big)\Big|_{s=0}&=
\iint\limits_{\R\times\R}\breve\psi_\e(x,y)\,dE_B(x)K\,dE_A(y)\\[.2cm]
&-\iint\limits_{\R\times\R}\breve\psi_\e(x,y)\,dE_A(x)K\,dE_A(y).
\end{align*}
By \rf{gf}, this is equal to 
\ba
&\int\limits_0^\be\iint\limits_{\R\times\R}
\Big(\big(S^*_t\psi_\e\big)*Q_t\Big)(x)e^{{\rm i}ty}\,dE_B(x)K\,dE_A(y)\,dt+\\[.2cm]
&\int\limits_0^\be\iint\limits_{\R\times\R}
\Big(\big(S^*_t\psi_\e\big)*Q_t\Big)(y)e^{{\rm i}tx}\,dE_B(x)K\,dE_A(y)\,dt-\\[.2cm]
&\int\limits_0^\be\iint\limits_{\R\times\R}
\Big(\big(S^*_t\psi_\e\big)*Q_t\Big)(x)e^{{\rm i}ty}\,dE_A(x)K\,dE_A(y)\,dt-\\[.2cm]
&\int\limits_0^\be\iint\limits_{\R\times\R}
\Big(\big(S^*_t\psi_\e\big)*Q_t\Big)(y)e^{{\rm i}tx}\,dE_A(x)K\,dE_A(y)\,dt.
\end{align*}
It is easy to see that
\ba
\int\limits_0^\be\iint\limits_{\R\times\R}
\Big(\big(S^*_t\psi_\e\big)*Q_t\Big)(x)e^{{\rm i}ty}\,dE_B(x)K\,dE_A(y)\,dt=
\int\limits_0^\be\Big(\big(S^*_t\psi_\e\big)*Q_t\Big)(B)K\exp({\rm i}tA)\,dt
\end{align*}
and similar equalities hold for the other three integrals. 


Thus
\ba
\psi_e(B)\!-\!\psi_\e(A)\!-\!\frac{d}{ds}\Big(\psi_\e(A_s)\Big)\Big|_{s=0}\!&=\!
\int\limits_0^\be\!\!
\Big(\!\big((S^*_t\psi_\e)*Q_t\big)(B)\!-\!\big((S^*_t\psi_\e)*Q_t\!\big)(A)\Big)
K\exp({\rm i}tA)dt\\[.2cm]
&+\int\limits_0^\be
\Big(\exp({\rm i}tB)-\exp({\rm i}tA)\Big)K\Big(\big((S^*_t\psi_\e)*Q_t\big)(A)\Big)\,dt.
\end{align*}

We have
$$
\lim_{\e\to0}\Big(\big(S^*_t\psi_\e\big)*Q_t\Big)(A)=\Big(\big(S^*_t\psi\big)*Q_t\Big)(A)
$$
in the strong operator topology (see the proof of Theorem 4 of \cite{Pe1}). By Lemma \ref{uL},
$$
\lim_{\e\to0}\Big(\big((S^*_t\psi_\e)*Q_t\big)(B)\!-\!\big((S^*_t\psi_\e)*Q_t\!\big)(A)\Big)=
\big((S^*_t\psi)*Q_t\big)(B)\!-\!\big((S^*_t\psi)*Q_t\!\big)(A)
$$
in $\bS_2$.

It follows that 
\ba
\lim_{\e\to0}&\trace\Big(\big(\big((S^*_t\psi_\e)*Q_t\big)(B)-\big((S^*_t\psi_\e)*Q_t\big)(A)\big)
K\exp({\rm i}tA)\Big)\\[.2cm]
=&\trace\Big(\big(\big((S^*_t\psi)*Q_t\big)(B)-\big((S^*_t\psi)*Q_t\big)(A)\big)
K\exp({\rm i}tA)\Big)
\end{align*}
and
\ba
\lim_{\e\to0}&\trace
\Big(\big(\exp({\rm i}tB)-\exp({\rm i}tA)\Big)K\Big(\big((S^*_t\psi_\e)*Q_t\big)(A)\big)\Big)\\[.2cm]
=&\trace
\Big(\big(\exp({\rm i}tB)-\exp({\rm i}tA)\Big)K\Big(\big((S^*_t\psi)*Q_t\big)(A)\big)\Big).
\end{align*}
Thus 
\ba
&\lim_{\e\to0}\trace\left(\psi_e(B)-\psi_\e(A)-\frac{d}{ds}\Big(\psi_\e(A_s)\Big)\Big|_{s=0}\right)\\[.2cm]
&=\int\limits_0^\be\trace
\Big(\big(\big((S^*_t\psi)*Q_t\big)(B)\!-\!\big((S^*_t\psi)*Q_t\!\big)(A)\big)
K\exp({\rm i}tA)\Big)\,dt\\[.2cm]
&+\int\limits_0^\be\trace\Big(
\big(\exp({\rm i}tB)-\exp({\rm i}tA)\big)K\big(\big((S^*_t\psi)*Q_t\big)(A)\big)\Big)\,dt\\[.2cm]
&=\trace\Big(\psi(B)-\psi(A)-\frac{d}{ds}\Big(\psi(A_s)\Big)\Big|_{s=0}\Big),
\end{align*}
which proves \rf{trf}. $\bl$

\medskip

{\bf Proof of Lemma \ref{uL}.} We have
$$
f_n(B)-f_n(A)=\iint\limits_{\R\times\R\setminus\D}\breve{f}_n(x,y)\,dE_B(x)K\,dE_A(y)
=\iint\limits_{\R\times\R\setminus\D}\breve{f}_n(x,y)\,d\E K(x,y),
$$
where $\E$ is the spectral measure on the space $\bS_2$ defined by
$\E(\d\times\s)T=E_B(\d)TE_A(\s)$, $\d,\,\s\subset\R$, $T\in\bS_2$,
and $\D\subset\R\times\R$ is the diagonal: $\D=\{(x,x):~x\in\R\}$.
Then
$$
\Big\|\big(f_n(B)-f_n(A)\big)-\big(f(B)-f(A)\big)\Big\|_{\bS_2}^2=
\iint\limits_{\R\times\R\setminus\D}\big|\breve{f}_n(x,y)-\breve{f}(x,y)\big|^2\,d(\E K,K)(x,y)\to0
$$
as $n\to\be$. $\bl$

Now we are going to to extend formula \rf{trf} to the whole class $B_{\be1}^2$. Consider first the case 
when $\f$ is a polynomial of degree at most 2. Clearly, for linear functions $\f$ the operator on the
left-hand side of \rf{trf} is the zero operator and the right-hand side of \rf{trf} is equal to 0.
Suppose now that $\f(t)=t^2$. If we perform formal manipulations, we obtain
\ba
&(A+K)(A+K)-A^2-\frac{d}{ds}(A+sK)(A+sK)\Big|_{s=0}\\[.2cm]
&=KA+AK+K^2-\frac{d}{ds}\Big(A^2+sKA+sAK+s^2K^2\Big)\Big|_{s=0}=K^2.
\end{align*}
We can put now by definition
$$
(A+K)^2-A^2-\frac{d}{ds}\Big(A+sK\Big)^2\Big|_{s=0}=K^2.
$$
The following result establishes formula \rf{trf} for the function $\f$, $\f(t)=t^2$.

\begin{thm}
\label{K2}
\bay
\label{eta}
\trace K^2=2\int_\R\eta(x)\,dx.
\ey
\end{thm}

\Pf To establish \rf{eta}, we first assume that $A$ is a bounded operator. Consider a sequence
$\{g_n\}_{n\ge1}$ such that 
$$
g_n(x)=x^2\quad\mbox{for}\quad x\in\big[-n,n\big],\quad\F g''_n\in L^1,
\quad\mbox{and}\quad\sup_{n\ge1}\|\F g''_n\|_{L^1}<\be.
$$
Then for $n\ge\|A\|+\|K\|$ we have
\ba
\trace K^2&=\trace\left(g_n(B)-g_n(A)-\frac{d}{ds}\Big(g_n(A_s)\Big)\Big|_{s=0}\right)\\[.2cm]
&=\int_\R g''_n(x)\eta(x)\,dx~
\to~2\int_\R\eta(x)\,dx,\quad\mbox{as}\quad n\to\be.
\end{align*}

If $A$ is an unbounded operator, consider the bounded self-adjoint operator $A_n$ defined by
$$
A_n=AE_A([-n,n]).
$$
Let $\eta_n$ be the generalized spectral shift function that correspond to the pair $(A_n, A_n+K)$.
Then
$$
\trace K^2=2\int_\R\eta_n(x)\,dx
$$
and \rf{eta} follows form the fact that
$$
\lim_{n\to\be}\int_\R\eta_n(x)\,dx=\int_\R\eta(x)\,dx,
$$
which can be found in \cite{Ko}. $\bl$

Finally, we obtain the following result.

\begin{thm}
\label{main}
The map 
$$
\f\mapsto\f(B)-\f(A)-\frac{d}{ds}\Big(\f(A_s)\Big)\Big|_{s=0}
$$
extends from $B^2_{\be1}(\R)\bigcap\Li$ to a bounded linear operator from $B^2_{\be1}(\R)$
to $\bS_1$ and trace formula {\em\rf{trf}} holds for functions $\f$ in $B^2_{\be1}(\R)$.
\end{thm}

\Pf Since the linear combinations of quadratic polynomials and functions whose Fourier
transforms have compact support in $\R\setminus\{0\}$ are dense in $B^2_{\be1}(\R)$,
the result follows immediately from Theorems \ref{int} and \ref{K2}. $\bl$

\

\section{\bf Open problems}
\setcounter{equation}{0}

\

The following interesting problems remain open.

\medskip 

{\bf Problem 1.} Suppose that $\f$ is a function of class $C^2$ on $\T$, $U$ is a unitary operator,
$A$ is a self-adjoint operator of class $\bS_2$. Is it true that
$$
\f(U_1)-\f(U)-\frac{d}{ds}\Big(\f(U_s)\Big)\Big|_{s=0}\in\bS_1,
$$
where $U_s=e^{{\rm i}sA}U$?

\medskip

{\bf Problem 2.} Suppose that $\f\in C^2(\R)$ and $\f''\in L^\be$. Let $A$ be a self-adjoint operator and
let $K$ be a self-adjoint operator of class $\bS_2$. Is it true that
$$
\f(A_1)-\f(A)-\frac{d}{ds}\Big(\f(A_s)\Big)\Big|_{s=0}\in\bS_1,
$$
where $A_s=A+sK$?

\medskip

Note that the right-hand sides of trace formulae \rf{ntf} and \rf{ktf} are well-defined for such
functions. I conjecture that the answer to both questions should be negative.

\

\

\noindent
\begin{tabular}{p{8cm}p{14cm}}
Department of Mathematics \\
Michigan State University  \\
East Lansing, Michigan 48824\\
USA
\end{tabular}


\begin{thebibliography}{99}
\bibitem[ABF]{ABF} {\sc J. Arazy, T. Barton}, and {\sc Y. Friedman}, 
{\em Operator differentiable functions} {\bf 13} (1990), 462--487.

\bibitem[BK]{BK} {\sc M.S. Birman} and {\sc M.G. Krein}, {\em On the theory of wave operators and
scattering operators}, Dokl. Akad. Nauk SSSR {\bf144} (1962), 475--478.

English transl.: Sov. Math. Dokl. {\bf3} (1962), 740--744.

\bibitem[BS1]{BS1} {\sc M.S. Birman} and {\sc M.Z. Solomyak}, {\em Double Stieltjes operator
integrals}, Problems of Math. Phys., Leningrad. Univ. {\bf1} (1966), 33--67 (Russian).

English transl.: Topics Math. Physics {\bf1} (1967), 25--54, Consultants Bureau Plenum
Publishing Corporation, New York.

\bibitem[BS2]{BS2} {\sc M.S. Birman} and {\sc M.Z. Solomyak}, {\em Double Stieltjes operator
integrals. II}, Problems of Math. Phys., Leningrad. Univ. {\bf2} (1967), 26--60 (Russian).

English transl.: Topics Math. Physics {\bf2} (1968), 19--46, Consultants Bureau Plenum
Publishing Corporation, New York.

\bibitem[BS3]{BS3} {\sc M.S. Birman} and {\sc M.Z. Solomyak}, {\em Remarks on the spectral
shift function},  Zapiski Nauchn. Semin. LOMI {\bf27} (1972), 33--46 (Russian).

English transl.: J. Soviet Math. {\bf3} (1975), 408--419.

\bibitem[BS4]{BS4} {\sc M.S. Birman} and {\sc M.Z. Solomyak}, {\em Double Stieltjes operator
integrals. III}, Problems of Math. Phys., Leningrad. Univ. {\bf6} (1973), 27--53 (Russian).

\bibitem[BS5]{BS5} {\sc M.S. Birman} and {\sc M.Z. Solomyak}, {\em Double operator integrals
in Hilbert space},  Int. Equat. Oper. Theory  {\bf47}  (2003), 131--168.

\bibitem[BY]{BY} {\sc M.S. Birman} and {\sc D.R. Yafaev}, {\em The spectral shift function. The
papers of M. G. Kre\u\i n and their further development},  Algebra i Analiz  {\bf4} 
(1992), 1--44 (Russian).  

English transl.: St. Petersburg Math. J.  {\bf4}  (1993), 833--870.

\bibitem[Bo]{Bo} J.-M. Bouclet, {\em Traces formulae for relatively Hilbert-Schmidt perturbations}, 
Asymptot. Anal.  {\bf32}  (2002), 257--291.

\bibitem[F]{F}  Yu.B. Farforovskaya, {\em An example of a Lipschitzian function of selfadjoint operators
that yields a nonnuclear increase under a nuclear perturbation}.  Zap. Nauchn. Sem. Leningrad. Otdel.
Mat. Inst. Steklov. (LOMI)  {\bf30}  (1972), 146--153 (Russian).

\bibitem[dPS]{dPS} B. de Pagter and F.A. Sukochev, {\em Differentiation of operator functions
in noncommutative \lb$L_p$-spaces}, Preprint, 2002.

\bibitem[Ko1]{Ko} {\sc L.S. Koplienko}, {\em The trace formula for perturbations of nonnuclear type},
Sibirsk. Mat. Zh. {\bf25:5}  (1984), 62--71 (Russian).

English transl.: Sib. Math. J. {\bf25} (1984), 735--743.

\bibitem[Ko2]{Ko2} {\sc L.S. Koplienko}, {\em Regularized spectral shift function for one-dimensional
Schr\"{o}dinger operator with slowly decreasing potential},   Sibirsk. Mat. Zh.  {\bf26:3}  (1985),  
72--77 (Russian).

English transl.: Sib. Math. J. {\bf26} (1985), 365--369.

\bibitem[Kr1]{Kr} {\sc M.G. Krein}, {\em On a trace formula in perturbation theory},
Mat. Sbornik {\bf33} (1953), 597--626 (Russian).

\bibitem[Kr2]{Kr2} {\sc M.G. Krein}, {\em On perturbation determinants and a trace formula for unitary and
self-adjoint operators},  Dokl. Akad. Nauk SSSR {\bf144} (1962) 268--271 (Russian).

\bibitem[Kr3]{Kr3} {\sc M.G. Krein}, {\em On some new investigations in the perturbation theory of
self-adjoint operators}, in: The First Summer Math. School, Kiev, 1964, 103--187 (Russian).

\bibitem[L]{L} {\sc I.M. Lifshitz}, {\em On a problem in perturbation theory 
connected with quantum statistics}, Uspekhi Mat. Nauk {\bf7} (1952), 171--180 (Russian).

\bibitem[N]{N} {\sc H. Neidhardt} {\em Spectral shift function and Hilbert--Schmidt
perturbation: extensions of some work of L.S. Koplienko}, Math. Nachr. {\bf 138} (1988), 7--25.

\bibitem[Pee]{Pe} {\sc J. Peetre}, 
{\em New thoughts on Besov spaces}, Duke Univ. Press., Durham, NC, 1976.

\bibitem[Pe1]{Pe1} {\sc V.V. Peller} {\em Hankel operators in the theory of perturbations
of unitary and self-adjoint operators},  Funktsional. Anal. i Prilozhen. {\bf19:2}  (1985), 
37--51 (Russian).

English transl.: Funct. Anal. Appl. {\bf19} (1985) , 111--123.

\bibitem[Pe2]{Pe2} {\sc V.V. Peller} {\em Hankel operators in the perturbation theory of 
of unbounded self-adjoint operators}.  Analysis and partial differential equations,  529--544,
Lecture Notes in Pure and Appl. Math., {\bf122}, Dekker, New York, 1990.

\bibitem[Pe3]{Pe3} {\sc V.V. Peller}, {\em When is a function of a Toeplitz operator close to a Toeplitz
operator?}  Toeplitz operators and spectral function theory,  59--85, Oper. Theory Adv. Appl., 42,
Birkh\"{a}user, Basel, 1989.

\bibitem[Pe4]{P} {\sc V.V. Peller}, {\em Hankel operators and their applications,}
Springer-Verlag, New York, 2003.

\bibitem[Y]{Y} {\sc D.R. Yafaev}, {\em Mathematical scattering theory. General theory}, Amer.
Math. Soc., Providence, 1992.  

\end{thebibliography}
\end{document}